\documentclass[a4paper,11pt,reqno,noindent]{amsart}
\usepackage[centertags]{amsmath}
\usepackage{amsfonts,amssymb,amsthm,dsfont,cases,amscd,esint,enumerate,fancyhdr,newlfont}
\usepackage[english]{babel}
\usepackage[body={15cm,21.5cm},centering]{geometry} 
\usepackage{tikz}
\usetikzlibrary{arrows.meta}

\theoremstyle{plain}
\newtheorem{theor0}{Theorem}
\newenvironment{theor}
  {\pushQED{\qed}\begin{theor0}}
  {\popQED\end{theor0}}

\numberwithin{equation}{section}

\newcommand{\e}{\varepsilon}
\newcommand{\dist}{\operatorname{dist}}
\newcommand{\R}{\mathbb R}
\newcommand{\Ic}{\mathcal I}
\newcommand{\Nc}{\mathcal N}
\newcommand{\cvf}{\rightharpoonup}
\newcommand{\loc}{{\operatorname{loc}}}
\newcommand{\Id}{\operatorname{Id}}
\newcommand{\E}{\mathbb{E}}

\newcommand{\D}{\operatorname{D}}
\newcommand{\Bb}{\bar{\boldsymbol B}}
\newcommand{\Ld}{\operatorname{L}}
\newcommand{\diam}{\operatorname{diam}}
\newcommand{\Div}{{\operatorname{div}}}
\newcommand{\Sym}{{\operatorname{sym}}}
\newcommand{\Skew}{{\operatorname{skew}}}
\newcommand{\Tr}{{\operatorname{tr}}}

\newcommand{\expec}[1]{\mathbb{E}\left[ #1 \right]}
\newcommand{\expecM}[1]{\mathbb{E}\bigg[ #1 \bigg]}

\title{Effective viscosity of semi-dilute suspensions}

\author[M. Duerinckx]{Mitia Duerinckx}
\address[Mitia Duerinckx]{Universit\'e Libre de Bruxelles, D\'epartement de Math\'ematique, 1050~Brussels, Belgium}
\email{mitia.duerinckx@ulb.be}
\author[A. Gloria]{Antoine Gloria}
\address[Antoine Gloria]{Sorbonne Universit\'e, CNRS, Universit\'e de Paris, Laboratoire Jacques-Louis Lions, 75005~Paris, France \& Institut Universitaire de France \& Universit\'e Libre de Bruxelles, D\'epartement de Math\'ematique, 1050~Brussels, Belgium}
\email{antoine.gloria@sorbonne-universite.fr}

\begin{document}

\begin{abstract}
This review is devoted to the large-scale rheology of suspensions of rigid particles in Stokes fluid. After describing recent results on the definition of the effective viscosity of such systems in the framework of homogenization theory, we turn to our new results on the asymptotic expansion of the effective viscosity in the dilute regime. This includes a new optimal proof of Einstein's viscosity formula for the first-order expansion, as well as the continuation of this expansion to higher orders. The essential difficulty originates in the long-range nature of hydrodynamic interactions: suitable renormalizations are needed and are captured by means of diagrammatic expansions.
\end{abstract}

\maketitle

\section{Introduction}
Suspensions of rigid particles in fluids are omnipresent in natural phenomena and in practical applications. They are known to display complex rheological behaviors on large scales, including possible non-Newtonian effects, e.g.~\cite{GM-11}, which we aim to understand and describe from a rigorous micro-macro perspective.
More precisely, we consider the macroscopic limit for a large number $N\gg1$ of small particles of size~$\e\ll1$ in a given tank~$U\subset\R^d$.
Neglecting both particle and fluid inertia, we assume that particles follow the fluid velocity and that the latter is instantaneously determined by the steady Stokes equations with no-slip conditions at particle boundaries.
Denoting by $\{I_{\e,N}^n\}_n$ the set of particles in $U$, with respective barycenters $\{x_{\e,N}^n\}_n$ and orientations $\{r_{\e,N}^n\}_n$, the steady Stokes equations for the fluid velocity $u_{\e,N}\in H^1_0(U)^d$ take the form
\begin{equation}\label{eq:Stokes}
-\triangle u_{\e,N}+\nabla p_{\e,N}=h,\qquad\Div(u_{\e,N})=0,\qquad\text{in $U\setminus\cup_nI_{\e,N}^n$},
\end{equation}
where $h\in\Ld^2(U)^d$ stands for some internal force in $U$. In view of no-slip conditions, we implicitly extend the fluid velocity $u_{\e,N}$  inside the particles, where it corresponds to the velocity of the particles. The rigidity of the particles then translates into
\begin{equation}\label{eq:Stokes-BC1}
\D(u_{\e,N})=0,\qquad\text{in $\cup_nI_{\e,N}^n$},
\end{equation}
where $\D(u)=\frac12(\nabla u+(\nabla u)')$ stands for the symmetric gradient. Equivalently, this means that for all $n$,
\begin{equation*}
u_{\e,N}^n=V_{\e,N}^n+W_{\e,N}^n(x-x_{\e,N}^n),\qquad\text{in $I_{\e,N}^n$},
\end{equation*}
for some translational velocity $V_{\e,N}^n\in\R^d$ and angular velocity tensor $W_{\e,N}^n\in\R^{d\times d}_\Skew$.
As we neglect the inertia of the particles, Newton's equations of motion reduce to the balance of forces and torques, which take the form of complementary boundary conditions,
\begin{eqnarray}\label{eq:Stokes-BC2}
\textstyle\e e+\fint_{\partial I_{\e,N}^n}\sigma(u_{\e,N},p_{\e,N})\nu&=&0,\nonumber\\
\textstyle\fint_{\partial I_{\e,N}^n}(x-x_{\e,N}^n)\times\sigma(u_{\e,N},p_{\e,N})\nu&=&0,\qquad\text{for all $n$},
\end{eqnarray}
where $\e e\in\R^d$ stands for some (weak) sedimentation force and \mbox{$\sigma(u,p):=2\D(u)-p\Id$} is the Cauchy stress tensor.
Given particle positions $\{I_{\e,N}^n\}_n$, the instantaneous fluid velocity~$u_{\e,N}$ is obtained as the unique solution of the above Stokes problem~\eqref{eq:Stokes}--\eqref{eq:Stokes-BC2}. Particle positions and orientations are then updated according to
\begin{equation*}
\partial_tx_{\e,N}^n=V_{\e,N}^n,\qquad\partial_t r_{\e,N}^n=W_{\e,N}^nr_{\e,N}^n,\qquad\text{for all $n$}.
\end{equation*}
In this way, particles interact via the fluid flow that they generate, and the resulting dynamics is reputedly complex in view of the multi-body, long-range, and singular nature of these hydrodynamic interactions.

Heuristically, while particles constitute small rigid inclusions in the fluid and hinder its flow, we expect homogenization to hold on large scales, leading to a notion of effective viscosity for the suspension. This effective viscosity naturally depends on the spatial arrangement of the particles, i.e.\@ on the microstructure, which evolves with the fluid flow and can thus adapt in time to external forces. This creates a possibly nonlinear response, hence non-Newtonian effects, which are indeed well-known in applications (e.g.\@ shear thinning of suspensions like ketchup). The mathematical understanding of such behaviors requires to couple homogenization with microstructure dynamics, which is still a fascinating open problem. Taking inspiration from the physics literature, we focus here on semi-dilute regimes and split the analysis into three steps:
\begin{enumerate}[$\bullet$]
\item {\it `Instantaneous' effective viscosity:}\\
Given particle positions $\{I_{\e,N}^n\}_n$, we expect the Stokes problem~\eqref{eq:Stokes}--\eqref{eq:Stokes-BC2} defining the fluid velocity to be approximated on large scales by an effective Stokes problem with some effective viscosity $\Bb$. This is by now well-understood in the framework of homogenization theory, and we review in Section~\ref{sec:eff-visc} our main results~\cite{DG-19,D-20a,DG-21a,DG-21b} on the topic. Keeping in mind the question of coupling homogenization with microstructure dynamics, we emphasize the importance of proving homogenization under the weakest possible assumptions on the microstructure.
\smallskip\item {\it Semi-dilute expansion of the effective viscosity:}\\
In the dilute regime, particles are sparse and interact little, hence the details of the microstructure should no longer be so relevant: we expect to expand the effective viscosity as $\Bb=\Id+\varphi\Bb^{(1)}+\ldots$ at low volume fraction $\varphi\ll1$, where the different terms would involve only reduced information on the microstructure that would be easier to track along the dynamics. To first order, this expansion is the celebrated {\it Einstein formula}, which has attracted considerable interest recently in the mathematical community.
In Section~\ref{sec:semi-dil}, we describe our new work~\cite{DG-22b} on the topic.
\smallskip\item {\it Coupling to semi-dilute microstructure dynamics:}\\
Coupling homogenization to microstructure dynamics remains an open problem even in the semi-dilute regime.
Only partial results are available for now, limited to mean-field regimes so dilute that they miss the description of non-Newtonian effects, cf.~\cite{Hofer-18,Mecherbet-19,Hofer-Schubert-21}. This is the subject of ongoing work that will not be discussed further in this review.
\end{enumerate}

\section{Effective viscosity problem}\label{sec:eff-visc}
We review recent results on the homogenization of the steady Stokes problem~\eqref{eq:Stokes}--\eqref{eq:Stokes-BC2} for given particle positions $\{I_{\e,N}^n\}_n$. More precisely, we shall consider a random ensemble for the set of particles, which is obtained by $\e$-rescaling of a given ensemble: given a random family $\{I_n\}_n$ of disjoint bounded subsets of $\R^d$, with respective barycenters $\{x_n\}_n$, we consider the following random set of rescaled inclusions in the reference domain $U$,
\begin{equation}\label{eq:Ieps-U}
\Ic_\e(U)\,:=\,\textstyle\bigcup_{n\in\Nc_\e(U)}\e I_n,\qquad\Nc_\e(U)\,:=\,\{n:\e I_n+\e B\subset U\},
\end{equation}
where $B$ stands for the unit ball, and we then consider the solution $u_\e\in H_0^1(U)^d$ of the corresponding Stokes problem
\begin{equation}\label{eq:ueps}
\left\{\begin{array}{ll}
-\triangle u_\e+\nabla p_\e=h,~\Div(u_\e)=0,&\text{in $U\setminus \Ic_\e(U)$},\\
\D(u_\e)=0,&\text{in $\Ic_\e(U)$},\\
|\e I_n|e+\int_{\e\partial I_n}\sigma(u_\e,p_\e)\nu=0,&\forall n\in\Nc_\e(U),\\
\int_{\e\partial I_n}(x-\e x_n)\times\sigma(u_\e,p_\e)\nu=0,&\forall n\in\Nc_\e(U).
\end{array}\right.
\end{equation}

\subsection{Structure of the equations}
Before moving on to the effective viscosity problem, we briefly comment on the Stokes problem~\eqref{eq:ueps} and emphasize that the different boundary conditions are natural for rigid particles: indeed, the weak formulation of~\eqref{eq:ueps} takes on the simple guise
\begin{multline*}
2\int_U\D(w):\D(u_\e)\,=\,\int_Uw\cdot\big(h\mathds1_{U\setminus\Ic_\e(U)}+e\mathds1_{\Ic_\e(U)}\big),\\
\text{for all $w\in H^1_0(U)^d$ with $\Div(w)=0$ and $\D(w)|_{\Ic_\e(U)}=0$}.
\end{multline*}
Here, $h\mathds1_{U\setminus\Ic_\e(U)}+e\mathds1_{\Ic_\e(U)}$ corresponds to the total force, combining the internal force in the fluid domain and the (weak) sedimentation force on the particles.
This weak formulation allows to represent
\begin{equation}\label{eq:represent-pi}
u_\e=\pi_\e v_\e,
\end{equation}
where $v_\e\in H^1_0(U)^d$ is the solution of the following Stokes problem without rigid particles,
\[-\triangle v_\e+\nabla q_\e=h\mathds1_{U\setminus\Ic_\e(U)}+e\mathds1_{\Ic_\e(U)},\qquad\Div(v_\e)=0,\qquad\text{in $U$},\]
and where $\pi_\e$ stands for the orthogonal projection 
\[\pi_\e~:~H^1_0(U)^d~\twoheadrightarrow~\big\{w\in H^1_0(U)^d:\Div(w)=0,\,\D(w)|_{\Ic_\e(U)}=0\big\},\]
with respect to the scalar product $(v,w)=\int_U \D(v):\D(w)$.
Another useful way to think of the above equations is to view rigid particles as droplets of diverging viscosity, cf.~e.g.~\cite{D-20a}: we have $u_\e^\theta\cvf{}u_\e$ in $H^1_0(U)^d$ as $\theta\uparrow\infty$, where $u_\e^\theta\in H^1_0(U)^d$ is the solution of the following Stokes problem with droplets of viscosity $\theta$,
\begin{equation}\label{eq:uepstheta}
-\Div\big(2(1+\theta\mathds1_{\Ic_\e(U)})\D(u_\e^\theta)\big)+\nabla p_\e^\theta=h\mathds1_{U\setminus\Ic_\e(U)}+e\mathds1_{\Ic_\e(U)},\quad\Div(u_\e^\theta)=0,\quad\text{in $U$}.
\end{equation}

\subsection{Homogenization}
In the macroscopic limit $\e\downarrow0$, provided that the ensemble of particles $\{I_n\}_n$ is stationary and ergodic, we expect the fluid velocity to be approximated by some effective Stokes problem, $u_\e\cvf{}\bar u$ in $H^1_0(U)^d$,
\begin{equation}\label{eq:homog-eqn}
-\Div(2\Bb\D(\bar u))+\nabla\bar p=(1-\varphi)h+\varphi e,\qquad\Div(\bar u)=0,\qquad\text{in $U$},
\end{equation}
where $\Bb$ is the effective viscosity tensor and where $(1-\varphi)h+\varphi e$ is the weak limit of the total force $h\mathds1_{U\setminus\Ic_\e(U)}+e\mathds1_{\Ic_\e(U)}$, in terms of the particle volume fraction~$\varphi=\expec{\mathds1_{\cup_nI_n}}$. The effective viscosity is naturally computed in terms of a cell problem: given a strain rate $E\in\R^{d\times d}_\Sym$ with $\Tr(E)=0$, we consider the correction $\Psi_E:=Ex+\psi_E$ of the linear straining flow $Ex$ that satisfies the following whole-space Stokes problem associated with the infinite family $\{I_n\}_n$ of rigid particles,
\begin{equation}\label{eq:corrector}
\left\{\begin{array}{ll}
-\triangle\Psi_E+\nabla\Sigma_E=0,~\Div(\Psi_E)=0,&\text{in $\R^d\setminus\cup_nI_n$},\\
\D(\Psi_E)=0,&\text{in $\cup_nI_n$},\\
\int_{\partial I_n}\sigma(\Psi_E,\Sigma_E)\nu=0,&\forall n,\\
\int_{\partial I_n}(x-x_n)\times\sigma(\Psi_E,\Sigma_E)\nu=0,&\forall n,
\end{array}\right.
\end{equation}
and the effective viscosity in direction $E$ is then given by the averaged dissipation rate
\begin{equation}\label{eq:eff-visc}
E:\Bb E\,=\,\expec{|\!\D(\psi_E)+E|^2}.
\end{equation}
As usual for cell problems in stochastic homogenization, the corrector $\psi_E$ solving the infinite-volume problem~\eqref{eq:corrector} is to be constructed in the class of random fields $\psi$ that are almost surely in $H^1_\loc(\R^d)^d$ such that the gradient $\nabla\psi$ is stationary, centered $\expec{\nabla\psi}=0$, and has finite second moments $\E[{|\nabla\psi|^2}]<\infty$. In particular, this entails that $\psi_E$ is a.s.\@ sublinear at infinity and is thus indeed a ``correction'' to the linear straining flow $Ex$.

Before our first work~\cite{DG-19} on this problem, the qualitative justification of this homogenization limit was still missing, surprisingly even in the periodic setting, as it required to deal at the same time with rigidity constraints and with incompressibility. Clearly, homogenization cannot hold without further geometric assumptions on the set of particles: if particles were to form an infinite chain, it would create macroscopic rigidity and~$\Bb$ would be infinite in the corresponding direction. The simplest way to prohibit such behaviors is the following, which has been largely used in previous work on the topic:
\begin{enumerate}[(H1)]
\item {\it Uniform condition:} assume $\rho_n:=\inf_{m:m\ne n}\dist(I_n,I_m)\ge\rho>0$ a.s.\@ for all $n$.
\end{enumerate}
This condition is quite restrictive and would not be preserved under microstructure dynamics, which motivates further generalizations. In~\cite{D-20a}, we showed that it can be replaced by a finer moment condition on interparticle distances:
\begin{enumerate}[(H2)]
\item {\it Moment condition:} assume $\E[\sum_n\rho_n^{-r_0}\mathds1_{I_n}]<\infty$ for some $r_0=r_0(d)$ large enough. (We can choose $r_0(d)=\frac32$ in dimension $d=3$, and $r_0(d)=0$ for $d>5$.)
\end{enumerate}
This condition is still not preserved under microstructure dynamics in 3D, cf.~\cite{BG-72a}. To weaken it further, the main difficulty is that we have poor control of the pressure in the vicinity of close particles, which hinders the standard proof of homogenization by Tartar's oscillating test function method, cf.~\cite{D-20a}.
Intuitively, however, the presence of close particles should not be problematic for the homogenization result: only long chains of close particles would be.
In~\cite{DG-21a}, we finally succeeded in avoiding the need for any detailed control at close particles, using instead a variational $\Gamma$-convergence approach. It holds under the following subcritical percolation type condition, which is essentially necessary for homogenization.
\begin{enumerate}[(H3)]
\item {\it Cluster condition:} Given some $\rho>0$, let $\{K_q\}_q$ be the family of connected components of the fattened set $\bigcup_nI_n+\rho B$, and assume that $\E[\sum_q(\diam K_q)^{r_0}\mathds1_{K_q}]<\infty$ for some $r_0=r_0(d)$ large enough.
\end{enumerate}
We can now state homogenization under any of the above conditions (H1)--(H3). Note that homogenization of stiff inclusions under~(H3) is new even in the scalar setting.

\begin{theor}[Qualitative homogenization; see~\cite{DG-19,D-20a,DG-21a}]\label{th:homog}
Assume that the ensemble of particles $\{I_n\}_n$ is stationary and ergodic, that particles are disjoint, uniformly bounded by \mbox{$I_n\subset B(x_n)$}, and uniformly $C^2$ almost surely. Further assume that either~\emph{(H1)}, \emph{(H2)}, or~\emph{(H3)} holds. Then the effective viscosity~\eqref{eq:eff-visc} is well-defined and we have $u_\e\cvf{}\bar u$ in~$H^1_0(U)^d$, where~$\bar u$ satisfies the effective Stokes problem~\eqref{eq:homog-eqn}.
\end{theor}

\subsection{Further results and extensions}
In~\cite{DG-21b}, under suitable mixing assumptions for the ensemble of particles $\{I_n\}_n$, we further prove quantitative error estimates for the homogenization limit, in form of
\[\textstyle\big\|u_\e-\bar u-\e\sum_{E}\psi_E(\tfrac\cdot\e)\partial_E\bar u\big\|_{\Ld^2(\Omega;H^1(U))}\,\lesssim_h\,\e^{\frac12}\times\left\{\begin{array}{lll}1&:&d>2,\\|\!\log\e|^{\frac14}&:&d=2,\end{array}\right.\]
where $\sum_{E}$ runs over an orthonormal basis of $\{E\in\R^{d\times d}_\Sym:\Tr(E)=0\}$.
We also obtain large-scale regularity results for the heterogeneous Stokes problem~\eqref{eq:ueps}, which play a key role in our work on sedimentation~\cite{DG-20}.
This builds on the quantitative homogenization theory recently developed for divergence-form uniformly elliptic problems, e.g.~\cite{AKM-book,GNO-reg}. A particular difficulty in our setting is that the field/flux structure of the Stokes problem is apparently destroyed due to rigid inclusions, in the sense that the na\"ive flux $\sigma(\Psi_E,\Sigma_E)\mathds1_{\R^d\setminus\cup_nI_n}$ is not divergence-free globally: instead, we take advantage of the representation~\eqref{eq:uepstheta} to recover the relevant notion of flux; see~\cite[Remark~4.2]{D-20a}.

We also mention~\cite{BDG-22} (see also~\cite{Girodroux-Lavigne-22}), where we generalize the above homogenization results to the case of suspensions of active swimmers. Swimming particles generate dipole forces in the fluid at small scales, which can lead to surprising rheological effects, such as a possible drastic reduction of the effective viscosity~\cite{Sokolov-Aranson-09}.

\section{Semi-dilute expansion}\label{sec:semi-dil}
While the effective viscosity $\Bb$ depends on the full microstructure $\{I_n\}_n$, cf.~\eqref{eq:eff-visc}, we expect this dependence to simplify perturbatively in the dilute regime. This idea takes its roots in the so-called effective medium expansions that emerged in the second half of the 19th century in the physics community; see e.g.~the historical account in~\cite{Markov-00}. For the effective viscosity problem, this was pioneered by Einstein~\cite{Einstein-05} in his PhD thesis, where he predicted that the first correction to the plain fluid viscosity is proportional to the particle volume fraction $\varphi$ and is given by the universal formula
\begin{equation}\label{eq:Einstein}
\Bb\,=\,\Id+\varphi\tfrac{d+2}2\Id+o(\varphi),\qquad\text{as $\varphi\downarrow0$},
\end{equation}
in case of spherical particles. Formally, this is obtained by summing single-particle contributions to the effective viscosity. The rigorous justification has attracted considerable interest in the mathematical community over the last decade, e.g.~\cite{Almog-Brenner,Haines-Mazzu,Niethammer-Schubert-19,Hillairet-Wu-19,GV-Hofer-20}, but two main questions have remained open:
\begin{enumerate}[---]
\item While previous contributions start from the (unphysical) uniform condition~(H1) on interparticle distances, can one also establish Einstein's formula~\eqref{eq:Einstein} under weaker assumptions of the form~(H2) or~(H3)?
\item What is the optimal error bound in~\eqref{eq:Einstein}? In particular, under what minimal assumption is the error indeed $o(\varphi)$?
\end{enumerate}
These questions are fully answered in our recent work~\cite{DG-22b}, as described in Section~\ref{sec:Einstein} below.
In link with the error bound in~\eqref{eq:Einstein}, there has also been interest in describing the next-order correction to Einstein's formula. Formally, the correction corresponds to summing pair contributions and the main difficulty is that this sum is not absolutely convergent: a renormalization is needed and a formal understanding was first achieved by Batchelor and Green~\cite{BG-72}. A few rigorous contributions have recently been devoted to this topic~\cite{GVH,GVM-20,GV-20}, although limited to special regimes.
The general description of higher-order corrections to Einstein's formula was first achieved in our recent work~\cite{DG-22b} in form of a cluster expansion, as described in Sections~\ref{sec:high}--\ref{sec:renorm} below. In the sequel, we focus on spherical particles $I_n=B(x_n)$ for simplicity, although it is not essential.

\subsection{Cluster expansion}\label{sec:cluster}
In the dilute regime, particles are typically well-separated and their interactions can thus formally be neglected perturbatively in the corrector problem~\eqref{eq:corrector}. The so-called cluster expansion corresponds to summing contributions of interactions between subsets of particles (or `clusters') of increasing cardinality.
For a finite index set $S\subset\R^d$, let $\psi_E^S$ be the solution of the corrector problem~\eqref{eq:corrector} associated with the finite set of particles $\{B(x)\}_{x\in S}$; in particular note that $\psi_E^\varnothing=0$. To first order, we then write
\[E:\Bb E\,\sim\,|E|^2+\text\guillemotleft~\expecM{\sum_n\big(|\!\D(\psi_E^{\{x_n\}})+E|^2-|E|^2\big)}~\text\guillemotright+\ldots\]
(We use quotation marks to indicate that this quantity is a priori not well defined, as indeed explained below.)
Using the short-hand notation
\[\delta^{\{x_n\}}|\!\D(\psi_E^{\#})+E|^2\,:=\,|\!\D(\psi_E^{\{x_n\}})+E|^2-|E|^2\]
for the first-order difference, and similarly defining higher-order differences, the formal cluster expansion takes the form
\begin{equation}\label{eq:cluster}
\Bb\,\sim\,\Id+\sum_{k=1}^\infty\Bb^{(k)},\qquad E:\Bb^{(k)}E\,=\,\text\guillemotleft~\expecM{\sum_{S\subset\{x_n\}_n\atop\sharp S=k}\delta^S|\!\D(\psi_E^{\#})+E|^2}~\text\guillemotright.
\end{equation}
The main difficulty is that for all $k\ge1$ the $k$th cluster term $\Bb^{(k)}$ is given by a series that is not absolutely convergent due to the long-range nature of hydrodynamic interactions. More precisely, to first order, we have
\begin{equation}\label{eq:lin-term-renorm}
\delta^{\{x_n\}}|\!\D(\psi_E^{\#})+E|^2\,=\,2E:\D(\psi_E^{\{x_n\}})+|\!\D(\psi_E^{\{x_n\}})|^2,
\end{equation}
where the single-particle solution satisfies $|\!\D(\psi_E^{\{x_n\}})|\simeq\langle \cdot-x_n\rangle^{-d}$, which entails that the definition of the first cluster term is indeed not absolutely convergent,
\begin{equation*}
\expecM{\sum_{n}\big|\delta^{\{x_n\}}|\!\D(\psi_E^{\#})+E|^2\big|}\,=\,\infty.
\end{equation*}
The same holds for all higher-order cluster formulas, cf.~\eqref{eq:cluster}, and suitable renormalization procedures are thus required.
Note, however, that divergence issues are only borderline and that cluster formulas can in fact be viewed as combinations of Calder\'on--Zygmund kernels --- this is explicit at first order in terms of the single-particle solution, cf.~\eqref{eq:lin-term-renorm}, but the structure is much more complicated and non-explicit in general.

\subsection{Correct scaling of cluster terms}\label{sec:scaling}
As the $k$th cluster term $\Bb^{(k)}$ involves contributions of $k$-tuples of particles, we might intuitively expect it to be of order $O(\varphi^k)$, but we argue that this cannot be true in general. By definition of the $k$-point density $f_k$ of the point process, the sum over $k$-tuples in the cluster formula~\eqref{eq:cluster} can formally be written as a multiple integral with respect to $f_k$,
\begin{equation*}
E:\Bb^{(k)}E\,=\,\text\guillemotleft~\int_{(\R^d)^k}\big(\delta^{\{x_1,\ldots,x_k\}}|\!\D(\psi_E^{\#})+E|^2(0)\big)\,f_k(x_1,\ldots,x_k)\,dx_1\ldots dx_k~\text\guillemotright,
\end{equation*}
or equivalently, by stationarity,
\begin{multline}\label{eq:cluster-fk}
E:\Bb^{(k)}E\,=\,\text\guillemotleft~\int_{(\R^d)^{k-1}}\Big(\int_{\R^d}\delta^{\{0,x_1,\ldots,x_k\}}|\!\D(\psi_E^{\#})+E|^2\Big)\\
\times\,f_k(0,x_1,\ldots,x_{k-1})\,dx_1\ldots dx_{k-1}~\text\guillemotright.
\end{multline}
Forgetting for now divergence issues at infinity,
and defining the $k$th-order intensity of the point process as
\begin{equation}\label{eq:lambda-k}
\lambda_k\,:=\,\|f_k\|_{\Ld^\infty((\R^d)^k)}
\end{equation}
(we refer to~\cite[Section~1.3.2]{DG-22b} for a more careful definition with local averages), we are led to expect, if the cluster formula makes any sense,
\[\Bb^{(k)}=O(\lambda_k).\]
In fact, due to long-range issues, the correct scaling is rather $\Bb^{(k)}=O(\lambda_k|\!\log\lambda_k|^{k-1})$ in general.
At first order, $\lambda_1=\lambda$ is the intensity of the point process and is of order $O(\varphi)$ as in Einstein's formula~\eqref{eq:Einstein}.
For a Poisson point process, due to tensorization, higher-order intensities are given by $\lambda_k=\lambda^k$, but this fails in general: for a strongly mixing point process, we only have
\[\lambda^j\le\lambda_j\le\lambda_{j-1}.\]
In particular, the correction to Einstein's formula is of order $O(\lambda_2|\!\log\lambda|)$ in general, so the approximation~\eqref{eq:Einstein} is only valid provided $\lambda_2|\!\log\lambda|=o(\lambda)$, which amounts to some weak form of local independence for the point process.

\subsection{Einstein's formula}\label{sec:Einstein}
With the above notation, we can now state our new main result on Einstein's formula. The novelty is the optimality of the error estimate, as well as the generality of the result, which holds under the weakest hypotheses for which homogenization is established, cf.~Theorem~\ref{th:homog}.

\begin{theor}[Einstein's formula; see Theorem~1 in~\cite{DG-22b}]\label{th:Einstein}
Under the same assumptions as in Theorem~\ref{th:homog}, we have
\[|\Bb-\Id-\Bb^{(1)}|\,\lesssim\,\lambda_2|\!\log\lambda|
+\left\{\begin{array}{ll}
0,&\text{under~\emph{(H1)}},\\
\lambda_2^{1-\frac1\kappa}\lambda^\frac1\kappa,&\text{under~\emph{(H2)} or~\emph{(H3)} with $\kappa=\kappa(r_0,d)$},
\end{array}\right.\]
where the first cluster term {\small$\Bb^{(1)}$} satisfies {\small$|\Bb^{(1)}|\simeq\lambda$} and is given by the well-defined renormalized cluster formula
\begin{equation}\label{eq:renorm-B1}
E:\Bb^{(1)}E\,:=\,\expecM{\sum_n|\!\D(\psi_E^{(x_n)})|^2}.
\end{equation}
In case of spherical particles, we recover Einstein's formula $\Bb^{(1)}=\varphi\frac{d+2}2\Id$.
\end{theor}

The renormalized formula~\eqref{eq:renorm-B1} is obtained from the cluster formula~\eqref{eq:cluster} by removing the linear term in~\eqref{eq:lin-term-renorm}. This is indeed natural:
starting from~\eqref{eq:cluster-fk} and expanding the difference, the cluster formula takes the form
\[E:\Bb^{(1)}E\,=\,\text\guillemotleft~\lambda\int_{\R^d}2E:\D(\psi_E^{\{0\}})~\text\guillemotright+\lambda\int_{\R^d}|\!\D(\psi_E^{\{0\}})|^2,\]
where the first linear term can only be meant to vanish as the integral of a gradient. As we shall see, this is to be understood for instance in a finite-volume approximation.

We now briefly describe our proof of the above result.
The idea is to neglect particle interactions, approximating the corrector locally by single-particle solutions, and we proceed by energy comparison. For that purpose, we denote by $\{V_n\}_n$ the Voronoi tessellation associated with the point process $\{x_n\}_n$, that is,
\[V_n\,:=\,\Big\{z\in\R^d:|z-x_n|<\inf_{m:m\ne n}|z-x_m|\Big\},\]
and for all $n$ we consider the solutions $\psi_{E;D}^{\{x_n\}}$ and $\psi_{E;N}^{\{x_n\}}$ of the Dirichlet and Neumann problems in $V_n$, respectively, associated with the single particle $I_n$. Comparing energies in Voronoi cells, we can bound the effective viscosity from above and below by infinite-volume averages of Dirichlet and Neumann single-particle energies, respectively,
\[\left\{\begin{array}{l}
E:\Bb E\,\le\,|E|^2+\lim_{R\uparrow\infty}|B_R|^{-1}\sum_{n:I_n\subset B_R}\int_{V_n}|\!\D(\psi_{E;D}^{\{n\}})|^2,\\
E:\Bb E\,\ge\,|E|^2+\lim_{R\uparrow\infty}|B_R|^{-1}\sum_{n:I_n\subset B_R}\int_{V_n}|\!\D(\psi_{E;N}^{\{n\}})|^2.
\end{array}\right.\]
As whole-space single-particle energies can also be bounded from above and below by Dirichlet and Neumann energies, we infer
\begin{multline*}
\Big|E:\Bb E-|E|^2-\lim_{R\uparrow\infty}|RB|^{-1}\sum_{n:I_n\subset RB}\int_{\R^d}|\!\D(\psi_E^{\{n\}})|^2\Big|\\
\,\le\,\lim_{R\uparrow\infty}|RB|^{-1}\sum_{n:I_n\subset RB}\Big(\int_{V_n}|\!\D(\psi_{E;D}^{\{n\}})|^2-\int_{V_n}|\!\D(\psi_{E;N}^{\{n\}})|^2\Big).
\end{multline*}
By elliptic regularity, the gap between Dirichlet and Neumann energies in $V_n$ can be controlled by $O(\rho_n^{-d})$, so that the above is bounded by $\E[\sum_n\rho_n^{-d}\mathds1_{I_n}]$. This can be evaluated in terms of the second-order intensity, and the conclusion follows.
Note that the generality of this proof makes it applicable to dilute systems in many other contexts; see also~\cite{DG-22c}.

\subsection{Higher-order cluster expansion}\label{sec:high}
The main difficulty to higher-order cluster expansions is that there is no simple way to make sense of cluster formulas~\eqref{eq:cluster} due to divergence issues. A natural idea is to start by considering finite-volume approximations of the effective viscosity,
\[E:\Bb_LE\,:=\,\expecM{\fint_{Q_L}|\!\D(\psi_{E;L})+E|^2},\]
where $\psi_{E;L}$ satisfies the corresponding corrector problem with periodic boundary conditions in the cube $Q_L:=(-\frac12L,\frac12L)^d$. The homogenization result, cf.~Theorem~\ref{th:homog}, yields in particular $\Bb_L\to\Bb$ in the infinite-volume limit $L\uparrow\infty$. For fixed $L$, as there is at most a finite number of particles in $Q_L$, cluster formulas make sense for all $k\ge1$,
\begin{equation}\label{eq:cluster-L}
E:\Bb_L^{(k)}E\,:=\,\expecM{\sum_{S\subset\{x_n\}_n\cap Q_L\atop\sharp S=k}\fint_{Q_L}\delta^S|\!\D(\psi_{E;L}^{\#})+E|^2},
\end{equation}
and we can investigate the associated cluster expansion.
It remains to prove good enough estimates on the different terms to pass to the infinite-volume limit.
In~\cite[Section~3]{DG-22b}, we prove two types of estimates,
\begin{enumerate}[$\bullet$]
\item {\it Direct estimates:} By a direct bound on cluster formulas~\eqref{eq:cluster-L}, we can capture the scaling in $\{\lambda_k\}_k$ as in Section~\ref{sec:scaling}, but the long-range decay of hydrodynamic interactions a priori leads to a logarithmic divergence in the infinite-volume limit,
\begin{equation}\label{eq:estim-direct}
|\Bb_L^{(k)}|\,\lesssim\,\lambda_k(\log L)^{k-1}.
\end{equation}
\item {\it Uniform-in-$L$ estimates:} Thinking of cluster formulas as complicated non-explicit combinations of Calder\'on--Zygmund kernels, we may expect to estimate them better by means of suitable energy estimates, carefully avoiding to take absolute values of the kernels. Taking inspiration from our previous work~\cite{DG-16a} on the effective conductivity problem, this can indeed be achieved by means of a hierarchy of interpolating energy estimates. Yet, in this way, we necessarily miss the scaling in $\{\lambda_k\}_k$,
\begin{equation}\label{eq:estim-unif}
|\Bb_L^{(k)}|\,\lesssim\,C^k.
\end{equation}
\end{enumerate}
As a direct consequence of uniform-in-$L$ estimates, we can deduce that infinite-volume limits
{\small$\Bb^{(k)}:=\lim_{L\uparrow\infty}\Bb_L^{(k)}$} exist, with convergence controlled by
\begin{equation}\label{eq:conv-BLk}
|\Bb_L^{(k)}-\Bb^{(k)}|\,\lesssim\,|\Bb_L-\Bb|^{2^{-k}}\,\to\,0,\qquad\text{as $L\uparrow\infty$}.
\end{equation}
This is shown in the same way as the fact that a sequence of functions that converge uniformly and have uniformly bounded derivatives also have converging derivatives.

It remains to interpolate between~\eqref{eq:estim-direct} and~\eqref{eq:estim-unif} to prove corresponding estimates in the infinite-volume limit. For this, we can think of the following shortcut: as the divergence in~\eqref{eq:estim-direct} is only logarithmic, it could be compensated by any algebraic rate in~\eqref{eq:conv-BLk}. Now, appealing to quantitative homogenization theory, e.g.~\cite{AKM-book}, under an algebraic $\alpha$-mixing condition, finite-volume approximations of the effective viscosity satisfy $|\Bb_L-\Bb|\lesssim L^{-\gamma}$ for some $\gamma>0$. Combined with~\eqref{eq:conv-BLk}, this yields
\begin{equation}\label{eq:conv-BLk-bis}
|\Bb_L^{(k)}-\Bb^{(k)}|\,\lesssim\,L^{- 2^{-k}\gamma},
\end{equation}
so that~\eqref{eq:estim-direct} implies
\[|\Bb^{(k)}|\,\lesssim\,|\Bb_L^{(k)}|+|\Bb_L^{(k)}-\Bb^{(k)}|\,\lesssim\,\lambda_k(\log L)^{k-1}+L^{-2^{-k}\gamma},\]
and we deduce by optimization $|\Bb^{(k)}|\lesssim\lambda_k|\!\log\lambda|^{k-1}$.
In this way, we get the following.

\begin{theor}[Higher-order cluster expansion; see Theorem~5 in~\cite{DG-22b}]\label{th:cluster}
Let the same assumptions hold as in Theorem~\ref{th:homog}, assume~\emph{(H1)} for simplicity, and assume that the ensemble of particles~$\{I_n\}_n$ satisfies an $\alpha$-mixing condition with algebraic rate. Then, for all $K\ge1$, the cluster expansion of the effective viscosity holds in form of
\begin{equation}\label{eq:cluster-exp-error}
\Big|\Bb-\Id-\sum_{k=1}^K\Bb^{(k)}\Big|\,\lesssim\,\sum_{k=K}^{2K}\lambda_{k+1}|\!\log\lambda|^{k},
\end{equation}
where cluster terms are defined by infinite-volume approximation~\eqref{eq:conv-BLk} and satisfy
\[|\Bb^{(k)}|\lesssim\lambda_k|\!\log\lambda|^{k-1}.\qedhere\]
\end{theor}

\subsection{Explicit renormalization}\label{sec:renorm}
Although Theorem~\ref{th:cluster} essentially solves the problem of higher-order cluster expansions, it has several drawbacks:
\begin{enumerate}[---]
\item Cluster terms are defined by infinite-volume approximation, cf.~\eqref{eq:conv-BLk}, which only provides an {\it implicit} renormalization of cluster formulas~\eqref{eq:cluster}. In particular, from this point of view, it is unclear whether logarithmic corrections in the estimates are actually optimal in the above statement.
\item The above relies on $\alpha$-mixing and quantitative homogenization theory as a black box, which seems however quite disconnected from the question of dilute expansions.
\item The convergence rate~\eqref{eq:conv-BLk-bis} deteriorates dramatically for large $k$, which is not expected to be optimal.
\end{enumerate}
All those criticisms led us to look for an explicit understanding of renormalization of cluster formulas, which was still open in the physics community beyond second order.

As explained after Theorem~\ref{th:Einstein}, the renormalization of the first-order cluster formula follows from the simple cancellation of the integral of a gradient in a finite-volume approximation, in form of {\small$\int_{Q_L}\D(\psi_{E;L}^{\{0\}})=0$}. At second order, the situation is already more complicated and we briefly recall the formal argument by Batchelor and Green~\cite{BG-72}: the cluster formula~\eqref{eq:cluster-fk} takes the form
\begin{equation*}
E:\Bb^{(2)}E\,=\,\text\guillemotleft~\int_{\R^d}\Big(\int_{\R^d}\delta^{\{0,y\}}|\!\D(\psi_E^{\#})+E|^2\Big)\,f_2(0,y)\,dy~\text\guillemotright,
\end{equation*}
or equivalently, after using corrector equations, cf.~\cite[Theorem~3]{DG-22b},
\begin{equation*}
E:\Bb^{(2)}E\,=\,\text\guillemotleft~\int_{\R^d}\Big(\int_{\partial B}\psi^{\{y\}}_E\cdot\sigma_E^{\{0,y\}})\nu\Big)\,f_2(0,y)\,dy~\text\guillemotright,
\end{equation*}
with the short-hand notation $\sigma_E^{S}:=\sigma(\psi_E^S+Ex,\Sigma_E^S)$. This can now be decomposed as
\begin{multline*}
E:\Bb^{(2)}E\,=\,\int_{\R^d}\Big(\int_{\partial B}\psi^{\{y\}}_E\cdot\big(\sigma_E^{\{0,y\}}-\sigma_E^{\{0\}}\big)\nu\Big)\,f_2(0,y)\,dy\\
+\text\guillemotleft~\int_{\R^d}\Big(\int_{\partial B}\psi^{\{y\}}_E\cdot\sigma_E^{\{0\}}\nu\Big)\,f_2(0,y)\,dy~\text\guillemotright,
\end{multline*}
where the first integrand has pointwise decay $|\int_{\partial B}\psi^{\{y\}}_E\cdot(\sigma_E^{\{0,y\}}-\sigma_E^{\{0\}})\nu|=O(\langle y\rangle^{-2d})$, so the first integral is absolutely convergent. Now noting that the second integrand has vanishing integral in a finite-volume approximation,
\begin{equation}\label{eq:cancel-Batchelor}
\int_{Q_L}\Big(\int_{\partial B}\psi_{E;L}^{\{y\}}\cdot\sigma_{E;L}^{\{0\}}\nu\Big)\,dy\,=\,\Big(\int_{Q_L}\psi_{E;L}^{\{y\}}\,dy\Big)\cdot\int_{\partial B}\sigma_{E;L}^{\{0\}}\nu\,=\,0,
\end{equation}
we can formally replace $f_2(0,y)$ by the two-point correlation \mbox{$h_2(0,y)=f_2(0,y)-\lambda^2$}, to the effect of
\begin{multline}\label{eq:Batchelor-Green}
E:\Bb^{(2)}E\,=\,\int_{\R^d}\Big(\int_{\partial B}\psi^{\{y\}}_E\cdot\big(\sigma_E^{\{0,y\}}-\sigma_E^{\{0\}}\big)\nu\Big)\,f_2(0,y)\,dy\\
+\int_{\R^d}\Big(\int_{\partial B}\psi^{\{y\}}_E\cdot\sigma_E^{\{0\}}\nu\Big)\,h_2(0,y)\,dy.
\end{multline}
The second integral is now absolutely convergent provided that $h_2$ has Dini decay (so that $y\mapsto\langle y\rangle^{-d}|h_2(0,y)|$ be integrable). This renormalized formula can be justified rigorously and leads to a fine understanding of the second cluster term.

\begin{theor}[Batchelor--Green renormalization; see Proposition~4.6 in~\cite{DG-22b}]
Let the same assumptions hold as in Theorem~\ref{th:homog}, assume~\emph{(H1)} for simplicity, and assume that the two-point correlation function satisfies $|h_2(0,y)|\le C\langle y\rangle^{-\gamma}$ for some $C,\gamma>0$. Then the second cluster term defined in~Theorem~\ref{th:cluster} is equivalently given by~\eqref{eq:Batchelor-Green}. Moreover, examination of this renormalized formula yields the following.
\begin{enumerate}[~$\bullet$]
\item The bound
{\small$|\Bb^{(2)}|\lesssim\lambda_2|\!\log\lambda|$}
is optimal in general, in the sense that it is attained by some point process. Yet, the logarithmic correction can be removed in some cases:
\begin{enumerate}[---]
\item if the point process is isotropic, we have ${\small|\Bb^{(2)}|}\lesssim\lambda_2$;
\item if the interparticle distance is $\ell:=\inf_{n\ne m}|x_n-x_m|\ge1$, we have ${\small|\Bb^{(2)}|}\lesssim\ell^{-2d}$.
\end{enumerate}
\smallskip\item The convergence rate~\eqref{eq:conv-BLk-bis} for finite-volume approximations can be improved to
\[|\Bb_L^{(2)}-\Bb^{(2)}|\,\lesssim\,L^{-\gamma}+\tfrac{\log L}L.\qedhere\]
\end{enumerate}
\end{theor}

The optimality of the logarithmic correction in the bound ${\small|\Bb^{(2)}|}\lesssim\lambda_2|\!\log\lambda|$ is understood as follows: the second term in~\eqref{eq:Batchelor-Green} can be written as $\int K(y)h_2(0,y)\,dy$ for some Calder\'on--Zygmund kernel $K$. As $\lambda_2$ is the sup norm of $h_2$, cf.~\eqref{eq:lambda-k}, the validity of the bound ${\small|\Bb^{(2)}|}\lesssim\lambda_2$ would require to control $K\ast$ in $\Ld^\infty(\R^d)$. This cannot be possible in general and we easily construct a correlation function $h_2$ for which it fails.

At higher orders, the renormalization of cluster formulas is more problematic as simple cancellations like~\eqref{eq:cancel-Batchelor} are no longer sufficient. Using corrector equations, the $k$th cluster formula~\eqref{eq:cluster-fk} can be written as
\begin{multline}\label{eq:cluster-fk-re}
E:\Bb^{(k)}E\,=\,\text\guillemotleft~\tfrac{k+1}2\int_{(\R^d)^{k-1}}\Big(\int_{\partial B}\delta^{\{x_1,\ldots,x_{k-1}\}}\psi_E^{\#}\cdot\sigma_E^{\{0\}}\nu\Big)\\
\times f_k(0,x_1,\ldots,x_{k-1})\,dx_1\ldots dx_{k-1}~\text\guillemotright+\ldots
\end{multline}
up to other similar terms that are slightly better behaved. In order to capture relevant cancellations, we introduce in~\cite{DG-22b} some new diagrammatic decomposition of corrector differences $\delta^{\{x_1,\ldots,x_{k-1}\}}\psi_E^{\#}$. Instead of using the classical method of reflections, which would yield an infinite series of terms only involving single-particle solutions, we introduce a suitable {\it finite} decomposition involving multiparticle solutions. For instance, we write
\[\nabla\delta^{x_1,x_2,x_3}\psi(0)
\,=\,
\begin{tikzpicture}[baseline={([yshift=-.8ex]current bounding box.center)},scale=0.4]
\begin{scope}[every node/.style={circle,fill,draw,inner sep=0pt,minimum size=3pt}]
    \node (2) at (0,0.6) {};
    \node (3) at (0,1.2) {};
    \node (4) at (0,1.8) {};
\end{scope}
\begin{scope}[every node/.style={circle,draw,fill=white,inner sep=0pt,minimum size=3pt}]
    \node (1) at (0,0) {};
\end{scope}
\begin{scope}[>={Stealth[black]},
every edge/.style={draw=black,thick}]
    \path [-] (1) edge (2);
    \path [-] (2) edge (3);
    \path [-] (3) edge (4);
\end{scope}
\end{tikzpicture}
+
\begin{tikzpicture}[baseline={([yshift=-.8ex]current bounding box.center)},scale=0.4]
\begin{scope}[every node/.style={circle,fill,draw,inner sep=0pt,minimum size=3pt}]
    \node (2) at (0,0.6) {};
    \node (3) at (0,1.2) {};
    \node (4) at (0,1.8) {};
\end{scope}
\begin{scope}[every node/.style={circle,draw,fill=white,inner sep=0pt,minimum size=3pt}]
    \node (1) at (0,0) {};
\end{scope}
\begin{scope}[>={Stealth[black]},
every edge/.style={draw=black,thick}]
    \path [-] (1) edge (2);
    \path [-] (2) edge (3);
    \path [-] (3) edge[bend left=30] (4);
    \path [-] (3) edge[bend left=-30] (4);
\end{scope}
\end{tikzpicture}
+
\begin{tikzpicture}[baseline={([yshift=-.8ex]current bounding box.center)},scale=0.4]
\begin{scope}[every node/.style={circle,fill,draw,inner sep=0pt,minimum size=3pt}]
    \node (2) at (0,0.6) {};
    \node (3) at (-0.4,1.2) {};
    \node (4) at (0.4,1.2) {};
\end{scope}
\begin{scope}[every node/.style={circle,draw,fill=white,inner sep=0pt,minimum size=3pt}]
    \node (1) at (0,0) {};
\end{scope}
\begin{scope}[>={Stealth[black]},
every edge/.style={draw=black,thick}]
    \path [-] (1) edge (2);
    \path [-] (2) edge[bend left=30] (3);
    \path [-] (2) edge[bend left=-30] (3);
    \path [-] (2) edge (4);
\end{scope}
\end{tikzpicture}
+
\begin{tikzpicture}[baseline={([yshift=-.8ex]current bounding box.center)},scale=0.4]
\begin{scope}[every node/.style={circle,fill,draw,inner sep=0pt,minimum size=3pt}]
    \node (2) at (0,0.6) {};
    \node (3) at (-0.4,1.2) {};
    \node (4) at (0.4,1.2) {};
\end{scope}
\begin{scope}[every node/.style={circle,draw,fill=white,inner sep=0pt,minimum size=3pt}]
    \node (1) at (0,0) {};
\end{scope}
\begin{scope}[>={Stealth[black]},
every edge/.style={draw=black,thick}]
    \path [-] (1) edge (2);
    \path [-] (2) edge[bend left=30] (3);
    \path [-] (2) edge[bend left=-30] (3);
    \path [-] (2) edge (4);
    \path [-] (4) edge (3);
\end{scope}
\end{tikzpicture}
+
\begin{tikzpicture}[baseline={([yshift=-.8ex]current bounding box.center)},scale=0.4]
\begin{scope}[every node/.style={circle,fill,draw,inner sep=0pt,minimum size=3pt}]
    \node (2) at (0,0.6) {};
    \node (3) at (-0.4,1.2) {};
    \node (4) at (0.4,1.2) {};
\end{scope}
\begin{scope}[every node/.style={circle,draw,fill=white,inner sep=0pt,minimum size=3pt}]
    \node (1) at (0,0) {};
\end{scope}
\begin{scope}[>={Stealth[black]},
every edge/.style={draw=black,thick}]
    \path [-] (1) edge (2);
    \path [-] (2) edge[bend left=30] (3);
    \path [-] (2) edge[bend left=-30] (3);
    \path [-] (2) edge[bend left=30] (4);
    \path [-] (2) edge[bend left=-30] (4);
\end{scope}
\end{tikzpicture}
+
\begin{tikzpicture}[baseline={([yshift=-.8ex]current bounding box.center)},scale=0.4]
\begin{scope}[every node/.style={circle,fill,draw,inner sep=0pt,minimum size=3pt}]
    \node (2) at (0,0.6) {};
    \node (3) at (-0.4,1.2) {};
    \node (4) at (0.4,1.2) {};
\end{scope}
\begin{scope}[every node/.style={circle,draw,fill=white,inner sep=0pt,minimum size=3pt}]
    \node (1) at (0,0) {};
\end{scope}
\begin{scope}[>={Stealth[black]},
every edge/.style={draw=black,thick}]
    \path [-] (1) edge (2);
    \path [-] (2) edge (3);
    \path [-] (3) edge (4);
    \path [-] (4) edge (2);
\end{scope}
\end{tikzpicture},\]
where the vertex
$\begin{tikzpicture}[baseline={([yshift=-.5ex]current bounding box.center)},scale=0.4]\begin{scope}[every node/.style={circle,draw,fill=white,inner sep=0pt,minimum size=3pt}]\node (1) at (0,0) {};\end{scope}\end{tikzpicture}$
stands for evaluation at $0$,
where other vertices correspond to integration variables $x_1,x_2,x_3$,
and where an edge between vertices $x$ and $y$ corresponds to a kernel with pointwise decay $\langle x-y\rangle^{-d}$.
More precisely, edges inside a given loop of the graph correspond to mutliparticle solutions (associated with the set of particles corresponding to the different vertices of the loop), while all the other edges correspond to single-particle solutions. This provides a crucial separation of variables.
These diagrammatic expansions are obtained by iterating Green's representation formula for corrector differences; see~\cite[Section~4.4.3]{DG-22b} for details.
We insert this decomposition into~\eqref{eq:cluster-fk-re} and we further expand the $k$-point density function~$f_k$ in terms of correlation functions, e.g.
\begin{equation*}
f_3(0,x_1,x_2)\,=\,\lambda^2+\lambda\big(h_2(0,x_1)+h_2(0,x_2)+h_2(x_1,x_2)\big)+h_3(0,x_1,x_2).
\end{equation*}
Assuming algebraic decay of correlations, this leads to additional couplings between the different integration variables, with some pointwise algebraic decay $\langle\cdot\rangle^{-\gamma}$, which we shall represent by dotted edges in the diagrams. Now the point is that our decomposition of corrector differences precisely allows to capture relevant cancellations: whenever a graph can be split into two subgraphs that are only connected by a simple edge, its contribution can be shown to vanish as a boundary term in the infinite-volume limit, e.g.
\[\begin{tikzpicture}[baseline={([yshift=-.8ex]current bounding box.center)},scale=0.4]
\begin{scope}[every node/.style={circle,fill,draw,inner sep=0pt,minimum size=3pt}]
    \node (2) at (0,0.6) {};
    \node (3) at (0,1.2) {};
    \node (4) at (0,1.8) {};
\end{scope}
\begin{scope}[every node/.style={circle,draw,fill=white,inner sep=0pt,minimum size=3pt}]
    \node (1) at (0,0) {};
\end{scope}
\begin{scope}[>={Stealth[black]},
every edge/.style={draw,thick}]
    \path [-] (1) edge (2);
    \path [-] (2) edge (3);
    \path [-] (3) edge (4);
\end{scope}
\begin{scope}[>={Stealth[black]},
every edge/.style={draw,thick,densely dotted}]
    \path [-] (1) edge[bend left=60] (3);
\end{scope}
\end{tikzpicture}~,
\begin{tikzpicture}[baseline={([yshift=-.8ex]current bounding box.center)},scale=0.4]
\begin{scope}[every node/.style={circle,fill,draw,inner sep=0pt,minimum size=3pt}]
    \node (2) at (0,0.6) {};
    \node (3) at (0,1.2) {};
    \node (4) at (0,1.8) {};
\end{scope}
\begin{scope}[every node/.style={circle,draw,fill=white,inner sep=0pt,minimum size=3pt}]
    \node (1) at (0,0) {};
\end{scope}
\begin{scope}[>={Stealth[black]},
every edge/.style={draw,thick}]
    \path [-] (1) edge (2);
    \path [-] (2) edge (3);
    \path [-] (3) edge (4);
\end{scope}
\begin{scope}[>={Stealth[black]},
every edge/.style={draw,thick,densely dotted}]
    \path [-] (1) edge[bend left=60] (2);
    \path [-] (3) edge[bend left=60] (4);
\end{scope}
\end{tikzpicture}~,
\begin{tikzpicture}[baseline={([yshift=-.8ex]current bounding box.center)},scale=0.4]
\begin{scope}[every node/.style={circle,fill,draw,inner sep=0pt,minimum size=3pt}]
    \node (2) at (0,0.6) {};
    \node (3) at (-0.4,1.2) {};
    \node (4) at (0.4,1.2) {};
\end{scope}
\begin{scope}[every node/.style={circle,draw,fill=white,inner sep=0pt,minimum size=3pt}]
    \node (1) at (0,0) {};
\end{scope}
\begin{scope}[>={Stealth[black]},
every edge/.style={draw,thick}]
    \path [-] (1) edge (2);
    \path [-] (2) edge[bend left=30] (3);
    \path [-] (2) edge[bend left=-30] (3);
    \path [-] (2) edge (4);
\end{scope}
\begin{scope}[>={Stealth[black]},
every edge/.style={draw,thick,densely dotted}]
    \path [-] (3) edge[bend left=-40] (1);
\end{scope}
\end{tikzpicture}
~=~0.\]
Starting from~\eqref{eq:cluster-fk-re} and removing all such contributions, we are left with terms that all correspond to absolutely convergent integrals, e.g.
\[E:\Bb^{(4)}E\,=\,
\begin{tikzpicture}[baseline={([yshift=-.8ex]current bounding box.center)},scale=0.4]
\begin{scope}[every node/.style={circle,fill,draw,inner sep=0pt,minimum size=3pt}]
    \node (2) at (0,0.6) {};
    \node (3) at (0,1.2) {};
    \node (4) at (0,1.8) {};
\end{scope}
\begin{scope}[every node/.style={circle,draw,fill=white,inner sep=0pt,minimum size=3pt}]
    \node (1) at (0,0) {};
\end{scope}
\begin{scope}[>={Stealth[black]},
every edge/.style={draw=black,thick}]
    \path [-] (1) edge (2);
    \path [-] (2) edge (3);
    \path [-] (3) edge (4);
\end{scope}
\begin{scope}[>={Stealth[black]},
every edge/.style={draw,thick,densely dotted}]
    \path [-] (1) edge[bend left=60] (3);
    \path [-] (2) edge[bend left=-60] (4);
\end{scope}
\end{tikzpicture}
+
\begin{tikzpicture}[baseline={([yshift=-.8ex]current bounding box.center)},scale=0.4]
\begin{scope}[every node/.style={circle,fill,draw,inner sep=0pt,minimum size=3pt}]
    \node (2) at (0,0.6) {};
    \node (3) at (0,1.2) {};
    \node (4) at (0,1.8) {};
\end{scope}
\begin{scope}[every node/.style={circle,draw,fill=white,inner sep=0pt,minimum size=3pt}]
    \node (1) at (0,0) {};
\end{scope}
\begin{scope}[>={Stealth[black]},
every edge/.style={draw=black,thick}]
    \path [-] (1) edge (2);
    \path [-] (2) edge (3);
    \path [-] (3) edge (4);
\end{scope}
\begin{scope}[>={Stealth[black]},
every edge/.style={draw,thick,densely dotted}]
    \path [-] (1) edge[bend left=60] (4);
\end{scope}
\end{tikzpicture}
+
\begin{tikzpicture}[baseline={([yshift=-.8ex]current bounding box.center)},scale=0.4]
\begin{scope}[every node/.style={circle,fill,draw,inner sep=0pt,minimum size=3pt}]
    \node (2) at (0,0.6) {};
    \node (3) at (0,1.2) {};
    \node (4) at (0,1.8) {};
\end{scope}
\begin{scope}[every node/.style={circle,draw,fill=white,inner sep=0pt,minimum size=3pt}]
    \node (1) at (0,0) {};
\end{scope}
\begin{scope}[>={Stealth[black]},
every edge/.style={draw=black,thick}]
    \path [-] (1) edge (2);
    \path [-] (2) edge (3);
    \path [-] (3) edge[bend left=30] (4);
    \path [-] (3) edge[bend left=-30] (4);
\end{scope}
\begin{scope}[>={Stealth[black]},
every edge/.style={draw,thick,densely dotted}]
    \path [-] (1) edge[bend left=60] (3);
\end{scope}
\end{tikzpicture}
+
\begin{tikzpicture}[baseline={([yshift=-.8ex]current bounding box.center)},scale=0.4]
\begin{scope}[every node/.style={circle,fill,draw,inner sep=0pt,minimum size=3pt}]
    \node (2) at (0,0.6) {};
    \node (3) at (0,1.2) {};
    \node (4) at (0,1.8) {};
\end{scope}
\begin{scope}[every node/.style={circle,draw,fill=white,inner sep=0pt,minimum size=3pt}]
    \node (1) at (0,0) {};
\end{scope}
\begin{scope}[>={Stealth[black]},
every edge/.style={draw=black,thick}]
    \path [-] (1) edge (2);
    \path [-] (2) edge (3);
    \path [-] (3) edge[bend left=30] (4);
    \path [-] (3) edge[bend left=-30] (4);
\end{scope}
\begin{scope}[>={Stealth[black]},
every edge/.style={draw,thick,densely dotted}]
    \path [-] (1) edge[bend left=60] (4);
\end{scope}
\end{tikzpicture}
+
\begin{tikzpicture}[baseline={([yshift=-.8ex]current bounding box.center)},scale=0.4]
\begin{scope}[every node/.style={circle,fill,draw,inner sep=0pt,minimum size=3pt}]
    \node (2) at (0,0.6) {};
    \node (3) at (-0.4,1.2) {};
    \node (4) at (0.4,1.2) {};
\end{scope}
\begin{scope}[every node/.style={circle,draw,fill=white,inner sep=0pt,minimum size=3pt}]
    \node (1) at (0,0) {};
\end{scope}
\begin{scope}[>={Stealth[black]},
every edge/.style={draw=black,thick}]
    \path [-] (1) edge (2);
    \path [-] (2) edge[bend left=30] (3);
    \path [-] (2) edge[bend left=-30] (3);
    \path [-] (2) edge (4);
\end{scope}
\begin{scope}[>={Stealth[black]},
every edge/.style={draw,thick,densely dotted}]
    \path [-] (1) edge[bend left=-50] (4);
\end{scope}
\end{tikzpicture}
+
\begin{tikzpicture}[baseline={([yshift=-.8ex]current bounding box.center)},scale=0.4]
\begin{scope}[every node/.style={circle,fill,draw,inner sep=0pt,minimum size=3pt}]
    \node (2) at (0,0.6) {};
    \node (3) at (-0.4,1.2) {};
    \node (4) at (0.4,1.2) {};
\end{scope}
\begin{scope}[every node/.style={circle,draw,fill=white,inner sep=0pt,minimum size=3pt}]
    \node (1) at (0,0) {};
\end{scope}
\begin{scope}[>={Stealth[black]},
every edge/.style={draw=black,thick}]
    \path [-] (1) edge (2);
    \path [-] (2) edge[bend left=30] (3);
    \path [-] (2) edge[bend left=-30] (3);
    \path [-] (2) edge (4);
\end{scope}
\begin{scope}[>={Stealth[black]},
every edge/.style={draw,thick,densely dotted}]
    \path [-] (1) edge[bend left=-50] (2);
    \path [-] (3) edge (4);
\end{scope}
\end{tikzpicture}
+
\begin{tikzpicture}[baseline={([yshift=-.8ex]current bounding box.center)},scale=0.4]
\begin{scope}[every node/.style={circle,fill,draw,inner sep=0pt,minimum size=3pt}]
    \node (2) at (0,0.6) {};
    \node (3) at (-0.4,1.2) {};
    \node (4) at (0.4,1.2) {};
\end{scope}
\begin{scope}[every node/.style={circle,draw,fill=white,inner sep=0pt,minimum size=3pt}]
    \node (1) at (0,0) {};
\end{scope}
\begin{scope}[>={Stealth[black]},
every edge/.style={draw=black,thick}]
    \path [-] (1) edge (2);
    \path [-] (2) edge[bend left=30] (3);
    \path [-] (2) edge[bend left=-30] (3);
    \path [-] (2) edge (4);
\end{scope}
\begin{scope}[>={Stealth[black]},
every edge/.style={draw,thick,densely dotted}]
    \path [-] (1) edge[bend left=50] (3);
    \path [-] (2) edge[bend left=-50] (4);
\end{scope}
\end{tikzpicture}
+
\begin{tikzpicture}[baseline={([yshift=-.8ex]current bounding box.center)},scale=0.4]
\begin{scope}[every node/.style={circle,fill,draw,inner sep=0pt,minimum size=3pt}]
    \node (2) at (0,0.6) {};
    \node (3) at (-0.4,1.2) {};
    \node (4) at (0.4,1.2) {};
\end{scope}
\begin{scope}[every node/.style={circle,draw,fill=white,inner sep=0pt,minimum size=3pt}]
    \node (1) at (0,0) {};
\end{scope}
\begin{scope}[>={Stealth[black]},
every edge/.style={draw=black,thick}]
    \path [-] (1) edge (2);
    \path [-] (2) edge[bend left=30] (3);
    \path [-] (2) edge[bend left=-30] (3);
    \path [-] (2) edge (4);
    \path [-] (4) edge (3);
\end{scope}
\begin{scope}[>={Stealth[black]},
every edge/.style={draw,thick,densely dotted}]
    \path [-] (1) edge[bend left=-50] (2);
\end{scope}
\end{tikzpicture}
+
\begin{tikzpicture}[baseline={([yshift=-.8ex]current bounding box.center)},scale=0.4]
\begin{scope}[every node/.style={circle,fill,draw,inner sep=0pt,minimum size=3pt}]
    \node (2) at (0,0.6) {};
    \node (3) at (-0.4,1.2) {};
    \node (4) at (0.4,1.2) {};
\end{scope}
\begin{scope}[every node/.style={circle,draw,fill=white,inner sep=0pt,minimum size=3pt}]
    \node (1) at (0,0) {};
\end{scope}
\begin{scope}[>={Stealth[black]},
every edge/.style={draw=black,thick}]
    \path [-] (1) edge (2);
    \path [-] (2) edge[bend left=30] (3);
    \path [-] (2) edge[bend left=-30] (3);
    \path [-] (2) edge (4);
    \path [-] (4) edge (3);
\end{scope}
\begin{scope}[>={Stealth[black]},
every edge/.style={draw,thick,densely dotted}]
    \path [-] (1) edge[bend left=50] (3);
\end{scope}
\end{tikzpicture}
+
\begin{tikzpicture}[baseline={([yshift=-.8ex]current bounding box.center)},scale=0.4]
\begin{scope}[every node/.style={circle,fill,draw,inner sep=0pt,minimum size=3pt}]
    \node (2) at (0,0.6) {};
    \node (3) at (-0.4,1.2) {};
    \node (4) at (0.4,1.2) {};
\end{scope}
\begin{scope}[every node/.style={circle,draw,fill=white,inner sep=0pt,minimum size=3pt}]
    \node (1) at (0,0) {};
\end{scope}
\begin{scope}[>={Stealth[black]},
every edge/.style={draw=black,thick}]
    \path [-] (1) edge (2);
    \path [-] (2) edge[bend left=30] (3);
    \path [-] (2) edge[bend left=-30] (3);
    \path [-] (2) edge (4);
    \path [-] (4) edge (3);
\end{scope}
\begin{scope}[>={Stealth[black]},
every edge/.style={draw,thick,densely dotted}]
    \path [-] (1) edge[bend left=-50] (4);
\end{scope}
\end{tikzpicture}
+\ldots\]
The number of terms in such decompositions grows very quickly, which probably makes them useless in practice, but it leads us to the proof of the following improved version of Theorem~\ref{th:cluster}, which is our main result in~\cite{DG-22b}.

\begin{theor}[Higher-order cluster expansion; see Proposition~4.8 in~\cite{DG-22b}]
Let the same assumptions hold as in Theorem~\ref{th:homog}, and assume~\emph{(H1)} for simplicity.
\begin{enumerate}[$\bullet$]
\item Given $K\ge1$, further assuming that correlation functions $h_{K+1},\ldots,h_{2K+1}$ have some decay $C\langle\cdot\rangle^{-\gamma}$ in each direction for some $C,\gamma>0$ (which is implied by algebraic $\alpha$-mixing), the cluster expansion of the effective viscosity holds in form of~\eqref{eq:cluster-exp-error},
where cluster terms are defined by infinite-volume approximation~\eqref{eq:conv-BLk}.
\smallskip\item For all $k\ge1$, further assuming that the correlation function $h_{k}$ has some decay $C\langle\cdot\rangle^{-\gamma}$ in each direction for some $C,\gamma>0$, then
\begin{enumerate}[---]
\item the $k$th cluster term is equivalently given by a renormalized formula only involving a finite number of absolutely convergent integrals;
\item it satisfies the bound $|\Bb^{(k)}|\lesssim\lambda_k|\!\log\lambda|^{k-1}$, which is optimal in general;
\item the convergence rate~\eqref{eq:conv-BLk-bis} for finite-volume approximations can be improved to
\[|\Bb_L^{(k)}-\Bb^{(k)}|\,\lesssim\,\tfrac{(\log L)^{k-1}}{L^{\gamma\wedge1}}.\qedhere\]
\end{enumerate}
\end{enumerate}
\end{theor}

\subsection{Remark: a new elliptic regularity result}
In order to prove the above, in~\cite{DG-22b}, we make repeated use of decay properties of solutions to Stokes problems with finite numbers of rigid particles. More precisely, we use for instance the following mean-value property, which seems new and might be of independent interest.

\begin{theor}[Mean-value property; see Appendix~A in~\cite{DG-22b}]
Let $\{I_n\}_{n\in S}$ be a finite collection of disjoint subsets of~$B_R$, assume that they are uniformly $C^2$ and that for some~\mbox{$\rho>0$} we have $\dist(I_n,I_m)\ge\rho$ and $\dist(I_n,\partial B_R)\ge\rho$ for all $n\ne m$.
Let $u\in H^1(B_R)^d$ satisfy the Stokes problem
\[\left\{\begin{array}{ll}
-\triangle u+\nabla p=0,~\Div(u)=0,&\text{in $B_R\setminus\cup_{n\in S}I_n$},\\
\D(u)=0,&\text{in $\cup_{n\in S}I_n$},\\
\int_{\partial I_n}\sigma(u,p)\nu=0,&\forall n\in S,\\
\int_{\partial I_n}(x-x_n)\times\sigma(u,p)\nu=0,&\forall n\in S.
\end{array}\right.\]
Then we have
\[\fint_B|\nabla u|^2\,\le\,C(\rho,\sharp S)\fint_{B_R}|\nabla u|^2,\]
where the multiplicative constant $C(\rho,\sharp S)$ only depends on $\rho$, $\sharp S$, and on dimension $d$.
\end{theor}

Note that this result is false in general if $\sharp S=\infty$, in link with classical counterexamples to Lipschitz regularity, but we show in~\cite{DG-21b} that a corresponding {\it annealed} estimate then holds upon taking a suitable ensemble average with respect to the set of particles.

\subsection{Further results and extensions}
In~\cite{BDG-22}, we derive similar results for the effective viscosity of {\it active} suspensions.
In a work in preparation with Pertinand, we further study Batchelor's dilute expansion~\cite{Batchelor-72} for the mean settling speed of a sedimenting suspension. The main difficulty is that, in the sedimentation problem, the long-range nature of hydrodynamic interactions is even more drastic: the flow disturbance at~$x$ due to a particle at~$y$ then decays like $\langle x-y\rangle^{1-d}$ instead of $\langle x-y\rangle^{-d}$. For this reason, a suitable renormalization is needed even to actually define the mean settling speed, which requires a nontrivial decay of correlations; see~\cite{DG-20}. The dilute expansion is then particularly tricky, but the same type of diagrammatic decompositions proves to be of crucial use in capturing cancellations.

\section*{Acknowledgements}
MD acknowledges financial support from F.R.S.-FNRS, and AG from the European Research Council (ERC) under the European Union's Horizon 2020 research and innovation programme (Grant Agreement n$^\circ$~864066).

\def\cprime{$'$} \def\cprime{$'$} \def\cprime{$'$}

\end{document}